\documentclass[a4paper,UKenglish,cleveref, autoref, thm-restate]{lipics-v2021}

\usepackage{todonotes}
\usepackage{xspace}
\usepackage{algcompatible}
\usepackage[noend]{algpseudocode}
\usepackage{algorithmicx, algorithm}
\usepackage{float}

\newcommand{\Vlabel}{\mathsf{Vlabel}}
\newcommand{\Tlabel}{\mathsf{Tlabel}}
\newcommand{\Over}{\mathsf{Over}}
\newcommand{\Under}{\mathsf{Under}}
\newcommand{\Full}{\mathsf{Full}}
\newcommand{\GLPartition}{\textsc{GL-Partition}\xspace}

\usepackage{nameref}

\newtheorem*{theorem*}{Theorem}

\makeatletter
\let\orgdescriptionlabel\descriptionlabel
\renewcommand*{\descriptionlabel}[1]{%
  \let\orglabel\label
  \let\label\@gobble
  \phantomsection
  \edef\@currentlabel{#1}%
  \let\label\orglabel
  \orgdescriptionlabel{#1}%
}
\makeatother

\nolinenumbers
\pdfoutput=1 
\hideLIPIcs  

\title{Connected Partitions via Connected Dominating Sets} 

\author{Aikaterini Niklanovits}{University Potsdam, Potsdam \and Hasso-Plattner-Institut, Potsdam}{Aikaterini.Niklanovits@hpi.de}{https://orcid.org/0000-0002-4911-4493}{}
\author{Kirill Simonov}{Department of Informatics, University of Bergen}{k.simonov@uib.no}{https://orcid.org/0000-0001-9436-7310}{}
\author{Shaily Verma}{Hasso-Plattner-Institut, Potsdam}{Shaily.Verma@hpi.de}{https://orcid.org/0009-0000-6789-1643}{}
\author{Ziena Zeif}{University Potsdam, Potsdam \and Hasso-Plattner-Institut, Potsdam}{Ziena.Zeif@hpi.de}{https://orcid.org/0000-0003-0378-1458}{}

\authorrunning{A. Niklanovits et al.} 

\Copyright{Jane Open Access and Joan R. Public} 



\category{} 

\relatedversion{} 

\EventEditors{John Q. Open and Joan R. Access}
\EventNoEds{2}
\EventLongTitle{42nd Conference on Very Important Topics (CVIT 2016)}
\EventShortTitle{CVIT 2016}
\EventAcronym{CVIT}
\EventYear{2016}
\EventDate{December 24--27, 2016}
\EventLocation{Little Whinging, United Kingdom}
\EventLogo{}
\SeriesVolume{42}
\ArticleNo{23}

\begin{document}

\maketitle
\begin{abstract}
The classical theorem due to Gy\H{o}ri and Lov\'{a}sz states that any $k$-connected graph $G$ admits a partition into $k$ connected subgraphs, where each subgraph has a prescribed size and contains a prescribed vertex, as long as the total size of target subgraphs is equal to the size of $G$. However, this result is notoriously evasive in terms of efficient constructions, and it is still unknown whether such a partition can be computed in polynomial time, even for $k = 5$.

We make progress towards an efficient constructive version of the Gy\H{o}ri--Lov\'{a}sz theorem by considering a natural strengthening of the $k$-connectivity requirement. 
Specifically, we show that the desired connected partition can be found in polynomial time, if $G$ contains $k$ disjoint connected dominating sets. As a consequence of this result, we give several efficient approximate and exact constructive versions of the original Gy\H{o}ri--Lov\'{a}sz theorem:
\begin{itemize}
    \item On general graphs, a Gy\H{o}ri--Lov\'{a}sz partition with $k$ parts can be computed in polynomial time when the input graph has connectivity $\Omega(k \cdot \log^2 n)$;
    \item On convex bipartite graphs, connectivity of $4k$ is sufficient;
    \item On biconvex graphs and interval graphs, connectivity of $k$ is sufficient, meaning that our algorithm gives a ``true'' constructive version of the theorem on these graph classes.
\end{itemize}
\end{abstract}
\section{Introduction}

 Partitioning a graph into connected subgraphs is a problem that naturally arises in various application areas such as image processing, road network decomposition and parallel graph systems for supporting computations on large-scale graphs \cite{DBLP:journals/dam/LucertiniPS93, DBLP:journals/jea/MohringSSWW06, DBLP:journals/tods/FanYXZLYLCX18, DBLP:conf/osdi/GonzalezLGBG12}.
Usually in such applications, there is a certain number of parts $k$ that the partition should have, as well as certain constraints on the size of each part.
In particular, a prominent objective is for the size of the parts to be as equal as possible and, therefore, for the partition to be as balanced as possible, known as Balanced Connected Partition problem ($BCP_k$).
This problem is proven to be NP-complete \cite{DBLP:journals/sigact/GareyJ78} and even hard to approximate \cite{CHLEBIKOVA1996225} and, hence, it has been studied in restricted graph classes such as grid graphs for which it was again proven to be NP-hard \cite{DBLP:journals/networks/BeckerLLS98}.

From a theoretical point of view, partitioning a graph into connected components each of a certain, chosen in advance, size is of independent interest as it gives us an insight into the structure of the graph.
With the simple example of a star graph, one can observe that it is not possible to partition it into connected sets where more than one of them has size $2$ or more.
In general, if there exists a set of $t$ vertices whose removal disconnects the graph (called a \emph{separator}), and if the number of parts $k$ is greater than $t$, then the $k$ parts of a connected partition might not necessarily be of arbitrary size.
Graphs that do not contain a separator of size strictly less than $t$ are called \emph{$t$-connected}.

The essence of this observation was captured in 1976 by Gy\H{o}ri \cite{Gyori1976} and Lov\'{a}sz \cite{Lovasz1977} who independently proved that the $k$-connectivity of a graph is a necessary and sufficient condition for any partition of a specific size to be possible, even when each part is additionally required to contain a specified vertex called a \emph{terminal}. This was later generalized to weighted graphs by Chandran et al. \cite{DBLP:conf/icalp/ChandranCI18}, Chen et al. \cite{DBLP:journals/jacm/ChenKLRSV07} and Hoyer \cite{Hoyer2016}.
Formally, we define the following computational problem.
\begin{definition}\label{def:GLPartition}
    In the \GLPartition problem, the input is a graph $G=(V, E)$, an integer $k \geq 1$, a set of $k$ distinct vertices $C=\{c_1,\ldots,c_k\}$ of $G$, and $k$ integers $N=\{n_1,\ldots,n_k\}$ such that $\sum_{i=1}^k n_i=|V|$, in short denoted as $(G,k,C,N)$. The task is to compute a partition of $V$ into $k$ disjoint subsets $V_1,\ldots, V_k$ such that $|V_i|=n_i$, $c_i\in V_i$, and the induced subgraph $G[V_i]$ on $V_i$ is connected, for all $i\in[k]$.
\end{definition}
We call the partition satisfying the property above a \emph{Gy\H{o}ri--Lov\'{a}sz partition}, or a \emph{GL-partition} for short. The result of Gy\H{o}ri and Lov\'{a}sz can then be stated as follows.
\begin{theorem*}[Gy\H{o}ri--Lov\'{a}sz Theorem]
    Given an instance of \GLPartition where the graph is $k$-connected, a GL-partition always exists.
\end{theorem*}

Although the proof of Gy\H{o}ri~\cite{Gyori1976} is constructive, the resulting algorithm runs in exponential time.
To this day, it remains open whether finding a GL-partition of a given $k$-connected graph could be done in polynomial time in general.
On the other hand, a few polynomial algorithms have been developed for special cases.
The original research direction is concerned with smaller values of $k$.
Specifically, in 1990, Suzuki et al. \cite{DBLP:journals/ipl/ZuzukiTN90, GLcase3} provided a polynomial algorithm for the cases $k=2$ and $k=3$.
Later, in 1993, Wada et al. \cite{DBLP:conf/wg/WadaK93} extended the result for $k=3$, and, in 1997, Nakano et al. \cite{DBLP:journals/ipl/NakanoRN97} showed a linear-time algorithm for the case where $k=4$ and the graph is planar with all the terminals located on the same face of the planar embedding.
More recently, in 2016, Hoyer and Thomas \cite{Hoyer2016} presented a polynomial-time algorithm for the case $k=4$ without additional restrictions on the graph structure. It is still open whether a polynomial-time algorithm exists for any fixed value of $k$ starting from $5$.
Another research direction focuses on polynomial-time constructions for the Gy\H{o}ri--Lov\'{a}sz theorem on specific graph classes, for all values of $k$.
Recently, Casel et al.~\cite{DBLP:conf/wg/CaselFINZ23} developed a polynomial-time algorithm for chordal graphs and (with a small additive violation of the size constraint) for $\text{HHI}^2_4$-free graphs. Additionally, Bornd\"{o}rfer et al.~\cite{DBLP:conf/approx/BorndorferCINSZ21} have shown an approximate version of the Gy\H{o}ri--Lov\'{a}sz theorem on general graphs, where the size of each part is within a factor of $3$ from the given target.

In the light of the challenges above, it is natural to ask whether efficient constructions exist for more restrictive versions of the Gy\H{o}ri--Lov\'{a}sz theorem.
Censor-Hillel, Ghaffari, and Kuhn~\cite{DBLP:conf/soda/Censor-HillelGK14} pioneered the notion of a connected dominating set (CDS) partition as a proxy for vertex connectivity of the graph. Formally, a graph admits a CDS partition of size $k$, if there exists a partition of the vertices into $k$ sets, such that each set forms a dominating set of the graph, and induces a connected subgraph.
It is easy to see that any graph that admits a CDS partition of size $k$, is also $k$-connected. The converse, unfortunately, does not hold; however, Censor-Hillel, Ghaffari, and Kuhn show that there is only a polylogarithmic gap between the connectivity and the size of the maximum CDS partition. Moreover, a suitable CDS partition can be efficiently computed.
The bound was later improved by Censor-Hillel et al.~\cite{10.1145/3086465}. 

\begin{theorem}[Corollary 1.6,~\cite{10.1145/3086465}]
    Every $k$-connected $n$-vertex graph has a CDS partition of size $\Omega(\frac{k}{\log^2 n})$, and such a partition can be computed in polynomial time.
    \label{thm:cds}
\end{theorem}
On the other hand, Censor-Hillel, Ghaffari, and Kuhn~\cite{DBLP:conf/soda/Censor-HillelGK14} show that a logarithmic gap between the connectivity and the size of a CDS partition is necessary.
\begin{theorem}[Theorem 1.3,~\cite{DBLP:conf/soda/Censor-HillelGK14}]
For any sufficiently large $n$, and any $k \ge 1$, there exist $k$-connected $n$-vertex graphs, where the maximum CDS partition size is $O(\frac{k}{\log n})$.
\label{thm:cds_lb}
\end{theorem}

More recently, Dragani\'{c} and Krivelevich~\cite{dragani2024disjoint} obtained a tighter bound of $(1 + o(1)) d/\ln d$ for CDS partitions in $d$-regular pseudorandom graphs.

\paragraph*{Our contribution}

The core result of our work is an efficient version of the Gy\H{o}ri--Lov\'{a}sz theorem that is based on the existence of a CDS partition of size $k$ in the given graph, which strengthens the original requirement of $k$-connectivity.

\begin{theorem}
There is an algorithm that, given an instance of \GLPartition together with a CDS partition of the graph of size $k$, computes a GL-partition in polynomial time.
\label{thm:cds_algo}
\end{theorem}

However, on its own Theorem~\ref{thm:cds_algo} would have a notable limitation, that computing an adequate CDS partition may in itself be an obstacle---as opposed to connectivity, which can always be computed exactly in polynomial time.
Fortunately, it turns out that a sufficiently large CDS partition can be found efficiently in a range of scenarios.
As the first example, Theorem~\ref{thm:cds} due to Censor-Hillel et al.~\cite{10.1145/3086465} provides such a construction for general graphs.
Combining it with our algorithm, we obtain the following efficient version of the Gy\H{o}ri--Lov\'{a}sz theorem, which requires higher connectivity, but puts no additional restrictions on the graph.

\begin{corollary}
\GLPartition is poly-time solvable on $\Omega(k \cdot \log^2 n)$-connected graphs.
\label{thm:general_algo}
\end{corollary}

Furthermore, we investigate the cases where Theorem~\ref{thm:cds_algo} can be applied to obtain a tighter connectivity gap.
Following the results of Censor-Hillel, Ghaffari, and Kuhn~\cite{DBLP:conf/soda/Censor-HillelGK14}, the case of general graphs cannot be improved beyond a logarithmic factor, as there exist $\Omega(k \cdot \log n)$-connected graphs that do not admit CDS partitions of size $k$, see Theorem~\ref{thm:cds_lb}. Therefore, putting additional restrictions on the input graph is necessary.

The first natural target here would be efficient constructions for the Gy\H{o}ri--Lov\'{a}sz theorem on certain graph classes in its exact version, that is, requiring connectivity exactly $k$. In the light of Theorem~\ref{thm:cds_algo}, it is sufficient to show a statement of the following form, that any $k$-connected graph in the class $\mathcal{G}$ admits a CDS partition of size $k$. This would imply that the Gy\H{o}ri--Lov\'{a}sz theorem admits an efficient construction on graphs in $\mathcal{G}$.
We identify two such cases, biconvex graphs and interval graphs.

\begin{restatable}{theorem}{cdsbiconvex}
    Let $G$ be a $k$-connected biconvex graph. A CDS partition of $G$ of size $k$ exists, and can be computed in polynomial time.
    \label{thm:cds_biconvex}
\end{restatable}
\begin{restatable}{theorem}{cdsinterval}
    Let $G$ be a $k$-connected interval graph. A CDS partition of $G$ of size $k$ exists, and can be computed in polynomial time.
    \label{thm:cds_interval}
\end{restatable}

Coupled with the algorithm of Theorem~\ref{thm:cds_algo}, this gives an efficient construction for the Gy\H{o}ri--Lov\'{a}sz theorem on biconvex graphs and interval graphs. 

\begin{corollary}
\GLPartition is poly-time solvable on $k$-connected biconvex graphs.
\label{thm:biconvex_algo}
\end{corollary}

We note that the efficient construction for chordal graphs was identified by Casel et al.~\cite{DBLP:conf/wg/CaselFINZ23}, implying the same result for interval graphs, since interval graphs form a subclass of chordal graphs. Nonetheless, we state the algorithmic consequence of Theorem~\ref{thm:cds_interval} next for completeness.

\begin{corollary}
\GLPartition is poly-time solvable on $k$-connected interval graphs.
\label{thm:interval_algo}
\end{corollary}

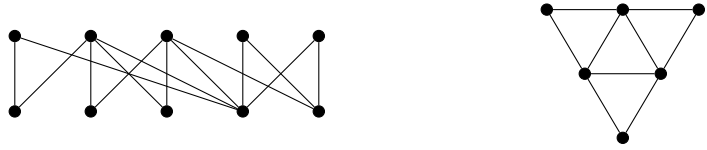
\begin{figure}
\begin{center}
    \begin{tikzpicture}[every node/.style={fill=black,inner sep=1.5pt}]
    \node[shape=circle,draw=black] (A) at (0,0) {};
    \node[shape=circle,draw=black] (B) at (1,0) {};
    \node[shape=circle,draw=black] (C) at (2,0) {};
    \node[shape=circle,draw=black] (D) at (3,0) {};
    \node[shape=circle,draw=black] (E) at (4,0) {};

    \node[shape=circle,draw=black] (A') at (0,-1) {};
    \node[shape=circle,draw=black] (B') at (1,-1) {};
    \node[shape=circle,draw=black] (C') at (2,-1) {};
    \node[shape=circle,draw=black] (D') at (3,-1) {};
    \node[shape=circle,draw=black] (E') at (4,-1) {};

    \path [-] (A) edge (A');
    \path [-] (A) edge (D');
    \path [-] (B) edge (A');
    \path [-] (B) edge (B');
    \path [-] (B) edge (C');
    \path [-] (B) edge (D'); 
    \path [-] (C) edge (B');
    \path [-] (C) edge (C');
    \path [-] (C) edge (D');
    \path [-] (C) edge (E');
    \path [-] (D) edge (D');
    \path [-] (D) edge (E');
    \path [-] (E) edge (D');
    \path [-] (E) edge (E');
    \begin{scope}[shift={(7,.35)}]
    \node[shape=circle,draw=black] (A) at (0,0) {};
    \node[shape=circle,draw=black] (B) at (1,0) {};
    \node[shape=circle,draw=black] (C) at (2,0) {};
    \node[shape=circle,draw=black] (D) at (0.5,-.85) {};
    \node[shape=circle,draw=black] (E) at (1.5,-.85) {};
    \node[shape=circle,draw=black] (F) at (1,-1.7) {};
    \end{scope}

    \path [-] (A) edge (B);
    \path [-] (A) edge (D);
    \path [-] (B) edge (C);
    \path [-] (B) edge (D);
    \path [-] (B) edge (E);
    \path [-] (C) edge (E); \path [-] (D) edge (E);
    \path [-] (D) edge (F);
    \path [-] (E) edge (F);
\end{tikzpicture}
\end{center}
\caption{Examples of a $2$-connected convex bipartite graph and a $2$-connected chordal graph that do not admit a CDS partition of size $2$}
\label{fig: counterexample}
\end{figure}

Unfortunately, there is a limit to the results of the form above, as already graph classes such as chordal graphs and convex bipartite graphs do not possess the property that $k$-connectivity implies the existence of a CDS partition of size $k$ (see Figure~\ref{fig: counterexample} for examples).
On the positive side, we are able to retain the productivity of Theorem~\ref{thm:cds_algo} even in such cases. We focus on the class of convex bipartite graphs and show that $4k$-connected convex bipartite graphs always admit a CDS partition of size $k$. 

\begin{restatable}{theorem}{cdsconvex}
    Let $G$ be a $4k$-connected convex bipartite graph. A CDS partition of $G$ of size $k$ exists, and can be computed in polynomial time.
    \label{thm:cds_convex}
\end{restatable}

As before, combined with the algorithm of Theorem~\ref{thm:cds_algo}, this gives us a ``$4$-approximate'' (in terms of connectivity) efficient version of the Gy\H{o}ri--Lov\'{a}sz theorem on this graph class, considerably improving the connectivity requirement of the general case.

\begin{corollary}
\GLPartition is solvable in polynomial time on $4k$-connected convex bipartite graphs.
\label{thm:convex_algo}
\end{corollary}


\paragraph*{Related work}
The problem of computing the largest CDS partition of a graph was first introduced by Hedetniemi and Laskar~\cite{HedetniemiL83}, where the maximum size of a CDS partition is referred to as the connected domatic number of the graph. A number of bounds on the size of a CDS partition are due to Zelinka~\cite{Zelinka1986}, including the fact that it is upper-bounded by the connectivity of the graph. Connected dominating sets have many practical applications for tasks such as routing and load-balancing in networks, for a detailed overview we refer to a survey due to Yu et al.~\cite{YuWWY13}.
Connections between CDS partitions and graph rigidity were explored by Krivelevich, Lew and Michaeli~\cite{KrivelevichLM23}.

Chordal graphs are graphs that do not admit induced cycles of length more than 3. The class of chordal graphs is well-studied, and is known to possess many useful properties; in particular, many problems that are NP-hard on general graphs turn out to be polynomial-time solvable on chordal graphs, such as computing chromatic number, independence number, treewidth~\cite{RoseTL76}. Interval graphs are intersection graphs of intervals on a real line, and form a subclass of chordal graphs~\cite{LekkeikerkerB62}.

Convex bipartite graphs\footnote{Often referred to simply as \emph{convex graphs} in the literature.} were introduced in 1967 by Glover~\cite{glover1967maximum} as they naturally arise in various industrial and scheduling applications. 
Formally, a bipartite graph $G=(A\cup B, E)$ is called \emph{convex} if there exists an ordering $\sigma$ of the vertices in $A$ such that the neighborhood of any vertex of $B$ is consecutive under $\sigma$.
Restricting to the class of convex bipartite graphs has been proven fruitful for many problems such as maximum independent set, minimum feedback set, (dynamic) maximum matching and independent dominating set~\cite{DBLP:journals/algorithmica/SoaresS09, DBLP:conf/mfcs/BrodalGHK07, DBLP:journals/acta/LiangC97, DBLP:journals/acta/LipskiP81, steiner1996linear, damaschke1990domination}.
On the other hand, several natural problems are known to be NP-complete on convex bipartite graphs. Most relevant to our work is the example of the $k$-path partition problem, proven to be NP-complete on convex bipartite graphs by Asdre and Nikolopoulos~\cite{ASDRE2007248}. Belmonte and Vatshelle~\cite{BelmonteV13} show that convex bipartite graphs have bounded mim-width, which provides a unified perspective on polynomial-time algorithms for various classes. Bonomo-Braberman et al.~\cite{BonomoBMP24} generalize this perspective further to extended notions of convexity.

Bipartite graphs with an ordering satisfying the convexity condition for both parts, are called biconvex. Biconvex graphs are thus a subclass of convex bipartite graphs. Abbas and Stewart studied properties of vertex orderings in biconvex graphs, yielding in particular that the vertex ranking problem is polynomial-time solvable on biconvex graphs~\cite{AbbasS00}. Apart from their structural and algorithmic interest, biconvex graphs and convex bipartite graphs have been proven useful in practice, as they model well situations such as task allocation problems~\cite{glover1967maximum, DBLP:journals/acta/LipskiP81}. In terms of other well-studied graph classes, convex bipartite graphs and biconvex graphs form a subclass of chordal bipartite graphs, which have a natural analogue of the chordal property among the bipartite graphs. For a more detailed introduction to the graph classes, their properties and relations between the classes, we refer to a survey by Brandst\"{a}dt, Bang Le and Spinrad~\cite{doi:10.1137/1.9780898719796}.

\subparagraph*{Paper organization.} In Section~\ref{sec:prelim}, we introduce preliminaries that we use throughout the paper. Section~\ref{section: All terminals on a tree} is dedicated to the main case of the algorithm constructing a Gy\H{o}ri--Lov\'{a}sz partition on the basis of a CDS partition, and Section~\ref{section: general case} completes the algorithm. We conclude with an overview and selected open problems in Section~\ref{sec:conclusion}. Due to space restrictions, the results for biconvex, interval, and convex bipartite graphs are deferred to appendix, as well as details from Section~\ref{section: general case} and selected proofs from Section~\ref{section: All terminals on a tree}. Such results are marked by ($\star$).


\section{Preliminaries}\label{sec:prelim}

All graphs mentioned in this paper are finite, undirected and simple.
We use standard graph theory notation, for reference see Appendix~\ref{sec:graph_notation} and the book of Diestel~\cite{DBLP:books/daglib/0030488}. We call a graph $G=(V,E)$ a \emph{path} if $V=\{v_1,\ldots,v_{n}\}$ and $E=\{v_iv_{i+1}|\text{for }i\in[n-1]\} $. We usually denote a path by $P$ and we call its order, the \emph{length} of $P$.
We call a graph $G=(V,E)$ a \emph{cycle} if $V=\{v_1,\ldots,v_n\}$ and $E=\{v_iv_{i+1}, v_1v_n|\text{for }i\in[n-1]\}$.
We call a graph $G=(V,E)$ a \emph{tree} if it does not contain any cycle as a subgraph.
Throughout this paper a vertex set of a tree is usually denoted by $T$, while its edge set is denoted by $E(T)$.
Given two graphs $H$, $G$ we say that $G$ is \emph{$H$-free} if $G$ does not contain $H$ as an induced subgraph.
Given two graphs $G$ and $T$, we say that $T$ is a \emph{spanning tree} of $G$ if $V(G)=V(T)$ and $G[V(T)]$ induces a tree.

Given two vertex sets $A,B\in V(G)$ we say that $A$ is \emph{dominating} to $B$ if every vertex in $B$ is adjacent to some vertex in $A$. 
Also, when we say that a vertex set is dominating to some graph $G$, we mean that it is dominating to its vertex set $V(G)$.
A vertex set $A$ is called a \emph{connected dominating set} (CDS) if and only if it is a dominating set and $G[A]$ is connected. 
Given a graph $G$, a collection of vertex sets $\mathcal V=\{V_1,\ldots, V_l\}$ such that $V_i\subseteq V(G)$ is called a \emph{connected dominating set partition} (CDS partition) if and only if the sets in $\mathcal V$ are pairwise disjoint and $\cup_{U\in\mathcal V}U=V(G)$, and each of them is a connected dominating set of $G$.
We call $|\mathcal V|$ the \emph{size} of the CDS partition $\mathcal V$.
Throughout the paper, we also use an equivalent notion of vertex-disjoint dominating trees. Given a CDS partition of a graph it is immediate to get a set of vertex-disjoint dominating trees of the same size by considering the spanning tree of each set in the CDS partition. Moreover, given $k$ disjoint dominating trees, one may obtain a CDS partition of size $k$ by simply adding the remaining vertices arbitrarily to those trees.
\section{From Simple CDS Partition to GL-partition}\label{section: All terminals on a tree}


Before we move on to the general algorithm of Theorem~\ref{thm:cds_algo} that computes a GL-partition from a CDS partition of size $k$, we cover in this section the main subroutine of this algorithm.
Specifically, we present an algorithm (Algorithm~\ref{case-alg:dominating trees}) that, given as input an instance of \GLPartition and $k$ dominating trees such that all terminals are contained in the same dominating tree, returns a partial GL-partition with $l$ sets, for some $l \in [k]$, such that the remaining vertices admit a CDS partition of size $k - l$. 
The general algorithm, described in the next section, then uses this result as a subprocedure.
We first describe the intuition behind the main steps of Algorithm~\ref{case-alg:dominating trees}, then provide a more detailed description of these steps, and finally define a formal procedure for each step and show its correctness.

{\bf Intuition of Algorithm~\ref{case-alg:dominating trees}:} 
Let $T_1$ be the dominating tree that contains $k$ terminal vertices. 
First we add all the vertices of $T_1$ to one of the $k$ vertex sets while maintaining the connectivity of each set. 
There exist $k-1$ dominating trees that have no terminal vertex.
Now, we map these $k$ vertex sets to one of these $k-1$ dominating trees. We aim to fulfill the vertex set to its demand while using the corresponding dominating tree.  

We use an iterative procedure that first assigns (but does not add) neighbors for each vertex set to these $k-1$ dominating trees, and then we add assigned vertices to each vertex set whose demand is not met by the assignment. 
Now, to fulfill the demand of each such vertex set, we use its corresponding dominating tree while preserving the connectivity of the vertex set.
At the end of this procedure, each part is either full or contains a whole dominating tree (apart from the one that was not assigned to any tree). 
The remaining vertices are arbitrarily added to the non-full sets, as each contains a dominating tree, apart from the one that has not been assigned a tree, for which we have ensured that there are enough adjacent vertices left.

\vspace{0.2cm}

{\bf Description of Algorithm~\ref{case-alg:dominating trees}:} 
The algorithm receives an instance of \GLPartition, $(G,k,C,N)$, where $C=\{c_1,\ldots,c_k\},\ N=\{n_1,\ldots,n_k\}$ and a set of dominating trees $\mathcal T=\{T_1,\ldots, T_k\}$ as an input such that such that $G$ is $k$-connected graph and $C\subseteq T_1$.

It is important to note that during this description of the algorithm, every time we add a vertex to a vertex set, it is implied that we also check whether the condition of having $l$ full sets say $V_1,\ldots, V_l$, intersecting exactly $l$ of the dominating trees, say $T_1,\ldots, T_l$ is met. If yes, then we output those $l$ full sets and return an instance of the original problem, $G-\cup_{i=1}^l V_i$ and $k-l$, $\mathcal T=\mathcal T\setminus\{T_{k-l},\ldots, T_k\}$, $C=C\setminus \{c_{k-l},\ldots, c_k\}$ and $N=N\setminus \{n_{k-l},\ldots, n_k\}$.
We moreover assume that every time we add a vertex we also check if a set meets its size demand and if so we add it to the set $\Full$ that includes all the connected vertex sets that satisfy their size demands. We denote this procedure of adding a vertex $u$ to set $V$ by $\text{Add}(V,u)$ (similarly $\text{Add}(V,U)$ when we have a vertex set $U$ instead of singleton $\{u\})$.
Similarly every time we remove a vertex from some set we check if it is in $\Full$ and if yes, we remove the set from $\Full$. We denote this procedure of removing vertex $u$ from the set $V$ by $\text{Remove}(V,u)$.


\begin{description}
    \item [Step 1\label{Case-Description: Add vertices}]
    We initialize each set $V_i$ to contain the corresponding terminal vertex $c_i$ for all $i\in [k]$ and then add the vertices of $T_1\setminus C$ one by one to a vertex set that it is adjacent to (based only on the adjacencies occurring by $T_1$).
    Next, we add the vertices of $V(G)\setminus V(\mathcal T)$ to an adjacent set one by one.
    
    \item [Step 2\label{Case-Description: Function Vlabel}]Each vertex of $\mathcal T\setminus T_1$ is assigned through a function $\Vlabel$ to an adjacent vertex set.
    Depending on whether the amount of vertices assigned to some vertex set is less than the vertices needed to be added for the set to satisfy its demand or not, the set is categorized in the set $\Under$ or $\Over$ respectively (Algorithm~\ref{case-alg:functions}).
    
    \item [Step 3\label{Case-Description: Function Tlabel}]Through the function $\Tlabel$ each set in $\Under$ is assigned one of the trees in $\mathcal T\setminus T_1$.
    Then, one vertex of each such tree adjacent to the corresponding set is added.
    Priority is given to the vertices assigned to this set through the function $\Vlabel$. However, if no such vertex exists, another adjacent vertex is arbitrarily chosen to be added, and we change the value of this vertex in the function $\Vlabel$ accordingly.
    If this results in a set moving from $\Over$ to $\Under$, we again assign a tree to this set (that has not been assigned to any other set) through the function $\Tlabel$ (Algorithm~\ref{case-alg:functions}).
    
    \item [Step 4\label{Case-Description: Add assigned trees}]For each tree assigned to some set, we consider its vertices not belonging in that set one by one and add them to that set as long as the set connectivity is preserved.
    Again, if, during this procedure, a set moves from $\Over$ to $\Under$, go to the previous Step.
    Repeat this procedure until each set in $\Under$ is either full or contains a whole tree (Algorithm~\ref{case-alg: Sets in Under contain a tree}).
    
    \item [Step 5\label{Case-Description: Add remaining vertices}]Add each of the remaining vertices, the ones not added to any set, to an adjacent non full set one by one (Algorithm~\ref{case-alg: Sets in Under contain a tree}).
    \item [Step 6\label{Case-Description: Restart}] Output the $k$ sets created.
\end{description}
As ~\ref{Case-Description: Add vertices} works in a straightforward way, we do not provide a pseudocode for it. We describe in the proof of the lemma below how this can be done in polynomial time. In the main algorithm, we refer to this procedure as $\text{AddTrees}(G,C,\mathcal T)$ that returns the sets $V_1,\ldots, V_k$.
\begin{restatable}[$\star$]{lemma}{lemmaAddTrees}\label{case-lemma: AddTrees}
     Given an instance of \GLPartition $(G,k,C,N)$ and $k$ vertex disjoint dominating trees $T_1,\ldots, T_k$ such that $C\subseteq T_1$, we construct $k$ vertex connected sets $V_1,\ldots, V_k$ such that $c_i\in V_i$ for all $i\in [k]$ and $T_1\subseteq\cup_{i=1}^k V_i$, in polynomial time.
\end{restatable}

Algorithm~\ref{case-alg:functions} realizes the two next steps of Algorithm~\ref{case-alg:dominating trees}. In particular, first each remaining vertex of the graph is assigned to some set through $\Vlabel$ and then each set is categorized in $\Over$ or $\Under$ depending on the number of vertices assigned to it.
Then, for each set in $\Under$ the assigned vertices are added and through $\Tlabel$ we assign a distinct tree each such set.
Afterwards we ensure that each set in $\Under$ contains a vertex from its assigned tree.

It is easy to notice that \ref{Case-Description: Function Vlabel} runs in polynomial time. 
In order to prove the correctness of \ref{Case-Description: Function Tlabel} in Lemma~\ref{case-lemma:ensure adjacencies in trees}, it suffices to show that there are at most $k-1$ sets in $\Under$ and this property is not affected while ensuring that each vertex set contains a vertex from the corresponding through $\Tlabel$ tree.
This part of the algorithm also runs in polynomial time because each set in $\Under$ is considered at most once and a set from $\Over$ might move to $\Under$ but once a set is in $\Under$ it will never move to $\Over$.
\begin{lemma}\label{case-lemma:ensure adjacencies in trees}
Let $(G,k,C,N)$ be an instance of \GLPartition, where $C=\{c_1,\ldots, c_k\}$, $\{V_1,\ldots, V_k\}$ are $k$ connected vertex sets such that $c_i\in V_i$ for all $i\in[k]$, and a set of $k$ vertex disjoint dominating trees $\mathcal T=\{T_1,\ldots, T_k\}$ such that $C\subseteq T_1$, are given as input to Algorithm~\ref{case-alg:functions}. Then, after the loop in line~\ref{Case-Step: Ensure adjacencies to assigned trees} (\ref{Case-Description: Function Tlabel}) each vertex set is either in $\Over$, or is assigned to a tree through $\Tlabel$ and contains a vertex from that tree.
    Moreover, $G[V_i\cup \Vlabel( V_i)]$ is connected for each $i\in [k]$ and at least one set is in $\Over$.
\end{lemma}
\begin{proof}
First notice that $G[V_i\cup \Vlabel( V_i)]$ is connected for all $i\in [k]$, by the definition of the function $\Vlabel$. We now show that there is at least one vertex set in $\Over$. 
Assume that there is no set in $\Over$ and hence $n_i-|V_i|>|\Vlabel(V_i)|$ for all $i\in [k]$. Then $\sum_{i=1}^kn_i-|V_i|>\sum_{i=1}^k|\Vlabel (V_i)|\Rightarrow n>\sum_{i=1}^k|\Vlabel (V_i)|+\sum_{i=1}^k|V_i|=n$ which is a contradiction.
This implies that regardless of whether a set moves from $\Over$ to $\Under$, as long as all vertices of $G$ are assigned to the sets $V_1,\ldots, V_k$, one of these sets will be in $\Over$.

In the loop in line~\ref{Case-Step: Ensure adjacencies to assigned trees}, we consider each tree assigned to some vertex set and first check whether this vertex set already contains some vertex of the tree. 
If yes, our lemma is satisfied. 
Otherwise, since the tree is dominating, there exists at least one vertex adjacent to the part of the vertex set intersecting $T_1$.
This vertex is already added to some vertex set in $\Under$ or assigned to one in $\Over$.
If it belongs in some set in $\Over$, then we simply add it to the vertex set.
Then we remove this vertex from the set that occurs from the assignment $\Vlabel$ for the set previously in $\Over$, and we check to see if it now meets the condition for this set belonging to $\Under$.
If yes, we go to line~\ref{Case-Step: Add vertices to Under} to deal with this set similarly. 
Repeating this procedure for the other sets in $\Under$ does not affect them as every vertex set currently in $\Under$ already contains all the vertices assigned to it through the function $\Vlabel$. 
Moreover, since in this step, we add a vertex to a vertex set, we also check if this set becomes full and intersects exactly one dominating tree. If that is true, we restart the algorithm for the adjusted values of $k$, $G$, $\mathcal T$, and $\mathcal C$.

\begin{algorithm}
\caption{{Labeling (\ref{Case-Description: Function Vlabel}, ~\ref{Case-Description: Function Tlabel})}\label{case-alg:functions}}

\begin{algorithmic}[1]
     \Require{An instance of \GLPartition $(G,k, C=\{c_1,\ldots, c_k\}, N=\{n_1,\ldots,n_k\})$, a set of $k$ vertex disjoint dominating trees $\mathcal T=\{T_1,\ldots, T_k\}$ such that $C\subseteq T_1$, and $k$ connected vertex sets $\{V_1,\ldots, V_k\}$ such that $c_i\in V_i$ for all $i\in [k]$}
    \State Initialize $\Vlabel:\{V_1,\ldots, V_k\}\rightarrow P(V(G))$ and $\Tlabel:\{V_1\ldots, V_k\}\to \mathcal T$ such that $\Vlabel(V_i)\gets\emptyset$, $\Tlabel(V_i)\gets \emptyset, \ \forall i\in[k]$.
    \State $\Over\gets\emptyset$, $\Under\gets\emptyset$
      \For{$i \gets 2 \textrm{ to } k$}
      \For {every vertex $v\in T_i$}\Comment{\ref{Case-Description: Function Vlabel}}
      \State Let $V_j\in \{V_1,\ldots, V_k\}$ such that $v\in N_G(V_j\cap T_1)$
      \State $\Vlabel(V_j)\gets \Vlabel(V_j)\cup v$    
      \EndFor
      \EndFor
      \For{$i \gets 2 \textrm{ to } k$}\label{Case-Step: Vertices Assignment}
      \If{$n_i-|V_i|\leq |\Vlabel(V_i)|$}
      \State $\Over\gets \Over\cup V_i$
      \Else
      \State $\Under\gets \Under\cup V_i$
      \EndIf
      \EndFor
      \For{every $V_i\in \Under$}\label{Case-Step: Add vertices to Under}\Comment{\ref{Case-Description: Function Tlabel}}
      \State $V_i\gets \text{Add}(V_i,\Vlabel(V_i))$
      \If{$\Tlabel(V_i)=\emptyset$}
      \State Let $T_j$ be a tree such that $T_j\not\in \cup_{s=1}^k\Tlabel (V_s)$
      \State $\Tlabel(V_i)\gets T_j$
      \EndIf
      \EndFor
      \For{every tree $T_i\in\mathcal T\setminus T_1$ such that $\Tlabel^{-1}(T_i)\neq\emptyset$}\label{Case-Step: Ensure adjacencies to assigned trees}
      \State Let $V_j$ be the vertex set such that $\Tlabel(V_j)=T_i$.
      \If{$V_j\cap T_i=\emptyset$}
      \State Let $u\in T_i$ such that $u\in N(V_j\cap T_1)$ and $V_s$ such that $u\in \Vlabel(V_s)$
      \If{$V_s\in \Over$}
      \State $V_j\gets \text{Add}(V_j,u)$, $\Vlabel(V_s)\gets \text{Remove}(\Vlabel(V_s), u)$.
      \If{$n_{s}-|V_s|>|\Vlabel (V_s)|$}
      \State $\Over\gets \Over \setminus V_s$, $\Under\gets \Under\cup V_s$
      \State Go to line ~\ref{Case-Step: Add vertices to Under}
      \EndIf
      \Else
      \State $V_s\gets \text{Remove}(V_s, u)$, $\Vlabel(V_s)\gets \Vlabel(V_s)\setminus u$, $V_j\gets \text{Add}(V_j, u)$
      \EndIf
      \EndIf
      \EndFor
      \State \Return $\Over, \Under, \Vlabel, \Tlabel$
    \end{algorithmic}
\end{algorithm}

Otherwise, if the adjacent vertex is part of some other set in $\Under$, we add it to our current set and remove it from the previous one (again checking if the condition for restarting the algorithm is met). 
Moreover, we remove it from the set created through the function $\Vlabel$ to avoid this vertex from existing in two sets (in case we go to line~\ref{Case-Step: Add vertices to Under}) in some other iteration of this loop.
Note that sets initially in $\Under$ might be added in $\Full$ or removed from it depending on whether we add or remove vertices from them, but they will never go from $\Under$ to $\Over$ as we need to maintain the property that each set in $\Over$ contains no vertices from $T_2,\ldots, T_k$.
\end{proof}

\begin{algorithm}
\caption{{Add Vertices (\ref{Case-Description: Add assigned trees},~\ref{Case-Description: Add remaining vertices})}
    \label{case-alg: Sets in Under contain a tree}}
    \begin{algorithmic}[1]
      \Require{An instance of \GLPartition $(G,k, C=\{c_1,\ldots, c_k\}, N=\{n_1,\ldots,n_k\})$, a set of $k$ vertex disjoint dominating trees $\mathcal T=\{T_1,\ldots, T_k\}$ such that $C\subseteq T_1$, $k$ connected vertex sets $\{V_1,\ldots, V_k\}$ such that $c_i\in V_i$ for all $i\in [k]$, sets $\Under,\ \Over,\ \Full$, and functions $\Vlabel,\ \Tlabel$.}
     \While{$\exists V_j\in \Under$ such that $V_j\not \in \Full$ and $T_i \not\subseteq V_j$ $\forall i \in \{2,\ldots,k\}$}\label{Case-Step: Sets in Under contain a tree}\Comment{\ref{Case-Description: Add assigned trees}}
      \State Let $T_i$ be the tree such that $\Tlabel (V_j)=T_i$
      \While{$V_j$ doesn't contain vertices only from $T_i$ and $T_i\not\subseteq V_j$ and $V_j\not\in \Full$}
      \State Let $v\in T_i$ such that $v\in N(V_j\cap(T_1\cup T_i))$
      \If{$v$ is assigned through the function $\Vlabel$ to some set $V_s\in \Over$}
      \State $V_j\gets \text{Add}(V_j, v)$, $\Vlabel (V_s)\gets \Vlabel(V_s)\setminus v$
      \If{$n_{s}-|V_s|>|\Vlabel (V_s)|$}
      \State $\Over\gets \Over \setminus V_s$, $\Under\gets \Under\cup V_s$
      \State Run Algorithm~\ref{case-alg:functions} with input $(G, \{T_1,\ldots, T_k\}, C, \Under, \Over, \Vlabel, N)$ starting from line 13.
      \EndIf
      \Else
      \State $V_s\gets \text{Remove}(V_s, v)$, $\Vlabel(V_s)\gets \Vlabel(V_s)\setminus v$, $V_j\gets \text{Add}(V_j,v)$
      \EndIf
      \EndWhile
      \EndWhile
        \For{every $V\in \Over$ and $V\not\in \Full$}\Comment{\ref{Case-Description: Add remaining vertices}}
        \State Let $u$ be a vertex in $\Vlabel(V)$ and not in $V$
        \State $V\gets \text{Add}(V, u)$
        \EndFor
        \For{every $v\in \cup_{i=2}^{k}T_i\setminus \cup_{i=1}^kV_i$}
        \State Let $V$ be a set such that $V\not\in \Full$ and $v\in N(V)$
        \State $V\gets \text{Add}(V, v)$
        \EndFor
        \State\Return $\{V_1,\ldots, V_k\}$
    \end{algorithmic}
\end{algorithm}

Algorithm~\ref{case-alg: Sets in Under contain a tree} realizes essentially the last two steps of the algorithm and either returns the desired partition of size $k$ or a smaller instance of a more general problem.
First we ensure that each set in $\Under$ is either full or contains the dominating tree assigned to it through $\Tlabel$. 
This is done by considering each set in $Under$ separately and then adding the neighboring vertices of the corresponding tree one by one. If during this procedure we "steal" a vertex assigned in some set in $\Over$ and this causes it to move to $\Under$ we run Algorithm~\ref{case-alg:functions} starting from line 13 (only running the loop of line 12 once, for the set that moved to $\Under$). This happens at most once for each set hence this part of Algorithm~\ref{case-alg: Sets in Under contain a tree} runs in polynomial time.
In the last part of this algorithm we simply add the remaining vertices to some adjacent non full set, first based on their assignment through $\Vlabel$ and then based on the adjacencies in $G$.

\begin{lemma}\label{case-lemma: Sets in Under contain a tree}
   Let an instance of \GLPartition $(G,k,C=\{c_1,\ldots,c_k\},N)$, $k$ connected vertex sets $\{V_1,\ldots, V_k\}$ such that $|V_i|<n_i$ and $c_i\in V_i$ for all $i\in[k]$, $k$ dominating to $G$ trees $\{T_1,\ldots, T_k\}$ and a function $\Vlabel: \{V_1,\ldots, V_k\} \to P(V(G))$, are given as an input to Algorithm~\ref{case-alg: Sets in Under contain a tree}. Then, after the loop in line~\ref{Case-Step: Sets in Under contain a tree} (\ref{Case-Description: Add remaining vertices}), each vertex set is connected and is either full or completely contains a dominating tree or is assigned through $\Vlabel$ to at least as many vertices as its remaining size demand that do not belong in any other vertex set.
\end{lemma}
\begin{proof}
    In order to see that the loop in line~\ref{Case-Step: Sets in Under contain a tree} terminates, it suffices to notice that it only considers sets in $\Under$ and the value of $\sum_{i=2}^{k}|T_i' \cap \Tlabel^{-1}(T_1)|$.

    We start by considering the vertex set $V_j$ such that $\Tlabel(V_j)=T_i$ for some value of $i$. 
    Due to Lemma~\ref{case-lemma:ensure adjacencies in trees}, we know that $T_i$ is either completely contained in $V_j$, in which case we do nothing and we proceed to consider the next tree, or it has a vertex of $T_i$ that is adjacent to $V_j\cap T_1$ and not in $V_j$.
    We then add that vertex to $V_j$ in the same manner as in the loop in line~\ref{Case-Step: Ensure adjacencies to assigned trees} of Algorithm~\ref{case-alg:functions} depending on whether this vertex belongs already in some set or not. 
    
    Notice that this action preserves the connectivity of the set $V_j$ due to the way the vertex is chosen and also does not affect the connectivity of other sets as this is based only on the adjacencies of the parts of each vertex set intersecting $T_1$, which is not affected.
    Moreover it increases the value of $\sum_{i=2}^{k}|T_i \cap \Tlabel^{-1}(T_1)|$ by $1$, which ensures that this loop terminates.
    Again, every time we add a vertex to a vertex set, it is implied that we also check if that set should meets its size demand and should be added to $\Full$ and also if the condition to terminate and output $l$ full sets and a smaller instance of the general problem is met.
    \end{proof}
Now that we have proved the correctness of \nameref{Case-Description: Add assigned trees}, notice that in the last part of Algorithm~\ref{case-alg: Sets in Under contain a tree} each non full set either contains a whole dominating tree or is adjacent to at least as many available vertices as it needs for its size demand to be met. 
We first start by the second case in which we add the assigned vertices to the set one by one until it is full. 
Then we are left with the case where each non full set contains a dominating tree, in which we can greedily add adjacent vertices to each such set one by one until all of them meet their size demands.

We now put everything together and provide the main algorithm of this section. 
\begin{algorithm}
\caption{{Gy\H{o}ri--Lov\'{a}sz theorem for special case of graphs with $k$ dominating trees}
    \label{case-alg:dominating trees}}
    \begin{algorithmic}[1]
       \Require{An instance of \GLPartition $(G,k, C, N)$ and a set of $k$ vertex disjoint dominating trees $\mathcal T=\{T_1,\ldots, T_k\}$ such that $C\subseteq T_1$.}
     \State $(V_1,\ldots, V_k)\gets \text{AddTrees}(G,C,\mathcal T)$
     \State $(\Over, \Under, \Vlabel, \Tlabel)\gets \text{Labeling}(G,C,\mathcal T, N, V_1,\ldots, V_k)$
     \State $(V_1,\ldots, V_k)\gets \text{Add Vertices}(G,C,\mathcal T, N, \Over, \Under,\Vlabel, \Tlabel, \Full)$
     \State Output $V_1,\ldots, V_k$ 
    \end{algorithmic}
    \end{algorithm}

\begin{restatable}[$\star$]{lemma}{lemmaCorrectness}\label{case-lemma: algorithm correctness}
    Let an instance of \GLPartition $(G,k,C=\{c_1,\ldots,c_k\},N=\{n_1,\ldots, n_k\})$ and $k$ vertex disjoint dominating trees $T_1,\ldots, T_k$ such that $C\subseteq T_1$ are given as an input to Algorithm~\ref{case-alg:dominating trees}.  Then Algorithm~\ref{case-alg:dominating trees} returns $l$ full sets $V_1,\ldots, V_l$ and a reduced instance, $(G-\cup_{i=1}^l V_i,k-l,C\setminus \{c_{k-l},\ldots, c_k\},N\setminus \{n_{k-l},\ldots, n_k\})$ of \GLPartition with a set of vertex-disjoint dominating trees of reduced graph $\mathcal T\setminus\{T_{k-l},\ldots, T_k\}$, for some $l\in[k]$, in polynomial time.

\end{restatable}

In the next section we tackle the more general instance in which there is a dominating tree that contains strictly less than $k$ terminals and some other contains at least one by showing that we are able to iteratively find a smaller instance in which all the terminals are contained in some dominating tree and then employ the algorithm developed in this section.

\section{From General CDS Partition to GL-partition}\label{section: general case}
\vspace{0.2cm}
In this section we provide the algorithm from Theorem~\ref{thm:cds_algo}.
The algorithm receives an instance of \GLPartition $(G,k,C,N)$ and $k$ vertex disjoint dominating trees $T_1,\ldots, T_k$ as input and it outputs a GL-partition $V_1,\ldots, V_k$ of $V(G)$.

If all terminal vertices are contained in one dominating tree then we can directly use Algorithm~\ref{case-alg:dominating trees} of the previous section.
Notice that if a terminal is not contained in any tree we can simply add it to one since each tree is dominating.
Hence, we assume without loss of generality that every terminal vertex is contained in some tree.
Moreover, we assume that the dominating trees are separated into three sets $\mathcal T_0,\mathcal T_1$ and $\mathcal T_{>1}$ depending on whether they contain $0$, $1$ or more than $1$ terminal vertices, respectively. This categorization can be done in polynomial time so for simplicity we assume that this is also part of the input.
As in the previous section we first provide an intuitive description of Algorithm~\ref{alg:dominating trees-general}.
The formal descriptions of the steps and the proof of correctness is deferred to the appendix.

{\bf Description of Algorithm~\ref{alg:dominating trees-general}:} 
The Algorithm receives as input a \GLPartition instance $(G,k,C,N)$, where $C=\{c_1,\ldots, c_k\}$, $N=\{n_1,\ldots, n_k\}$, and a set of $k$ vertex-disjoint dominating trees $\mathcal T=\{T_1,\ldots, T_k\}$ separated into three sets $\mathcal T_0,\mathcal T_1,\mathcal T_{>1}$ depending on how many terminal vertices they contain.
It is important to note that, like in the previous section, during this description of the algorithm, every time we add a vertex to a vertex set, it is implied that we also check if the condition of having $l$ full sets say $V_1,\ldots, V_l$, intersecting exactly $l$ of the dominating trees, say $T_1,\ldots, T_l$ is met. If yes, then we output those $l$ full sets and restart the algorithm for $G-\cup_{i=1}^l V_i$ and $k-l$, $\mathcal T=\mathcal T\setminus\{T_{k-l},\ldots, T_k\}$, $C=\{c_{k-l},\ldots, c_l\}$ and $N=\{n_{k-l},\ldots, n_k\}$. 
For simplicity and space constraints, we do not include this in every step of this description nor in the pseudocodes provided in this section.
We moreover assume both in pseudocodes of this
section and in this description that every time we add a vertex we also check if a set meets
its size demand and if so we add it to the set $\Full$ that includes all the connected vertex sets
that satisfy their size demands. Similarly every time we remove a vertex from some set we
check if it is in $\Full$ and if yes, we remove it from $\Full$.
\begin{description}
    \item [Step 1\label{Description: Add vertices}]We first initialize each vertex set $V_i$ to contain the corresponding terminal vertex $c_i$, where $i\in [k]$.
    Then, for each tree in $\mathcal T_1\cup \mathcal T_{>1}$, we add its vertices one by one to a vertex set that it is adjacent to (based only on the adjacencies occurring only by the tree considered).
    Next, we add the vertices of $V(G)\setminus V(\mathcal T)$ to an adjacent set one by one.
    \item [Step 2 \label{Description: Chose a tree set}] For each tree $T\in \mathcal T_{>1}$ that contains $l$ terminals, we create a tree set $\mathcal T'$ containing $T$ and $l-1$ trees from $\mathcal T_0$, which we then remove from $\mathcal T_{>1}$ and $\mathcal T_0$ respectively. 
    Among those sets, we choose one that has enough vertices available to satisfy the size demands of all of the $l$ vertex sets intersecting $T$.
    (Algorithm~\ref{alg:Separate trees to sets})
    \item [Step 3\label{Description: Case Algorithm}] Call Algorithm~\ref{case-alg:dominating trees} that has as input those $l$ dominating trees, the corresponding vertex sets and terminal vertices, while the graph we work on is the one that occurs from the union of those sets and trees. 
    Notice that the sum of the demands of those sets might be less than the size of the graph but this does not affect the correctness of the algorithm. 
\end{description}

\section{Conclusion}\label{sec:conclusion}


As we showcase in this work, the algorithm of Theorem~\ref{thm:cds_algo} opens up a new avenue for efficiently computing connected partitions, in the settings where sufficiently large CDS partitions are available.
This motivates further research into improving bounds for CDS partitions beyond the worst-case instances, as well as into parameterized and approximation algorithms for computing a CDS partition of the largest size.

To the best of our knowledge, Corollaries~\ref{thm:general_algo} and~\ref{thm:convex_algo} are the first examples of a Gy\H{o}ri--Lov\'{a}sz type of results that allow for a trade-off between efficiency of the construction and the connectivity requirement. It would be interesting to investigate whether a more refined trade-off is possible, interpolating smoothly between the polynomial time for higher connectivity and exponential time for extremal connectivity.

Our work also initiates the study of the gap between connectivity and the size of a CDS packing in structured graph classes.
It remains a curious open question to characterize the class of graphs, for which $k$-connectivity implies the existence of a $k$-CDS partition. So far, our results show that this class contains interval and biconvex graphs (Theorems~\ref{thm:cds_interval} and~\ref{thm:cds_biconvex}), but does not contain chordal and convex bipartite graphs.



\bibliography{bibliography}
\clearpage
\appendix
\section{Extended Preliminaries}\label{sec:graph_notation}

A graph $G=(V,E)$ is defined by its set of vertices $V$, also denoted by $V(G)$, and its set of edges $E$, also denoted by $E(G)$.
We say that $|V|$ is the \emph{order} of $G$, and we usually denote it by $n$.
Given a graph $G=(V,E)$, we denote by $P(V)$ the powerset of the vertices of $G$.
Given a vertex set $X\in V$ we denote by $N_G(X)=\{u\in V\setminus X:uv\in E \text{ and } v\in X\}$ the \emph{open neighbourhood} of $X$ in $G$ and by $N_G[X]=X\cup N(X)$ the \emph{closed neighbourhood} of $X$ in $G$.
If the graph we refer to, is clear by the context we omit $G$ as a subscript.
We also refer to $|N(v)|$ as the \emph{degree of the vertex $v$}, denoted $\deg(v)$.
Moreover, given a vertex set $X\in V(G)$ we denote by $G[X]$ the \emph{subgraph of $G$ induced by $X$}, the graph with vertex set $X$ and with edge set the subset of $E$ that contains all the edges with both endpoints in $X$.

We say that a graph is $k$-connected if the removal of any vertex set of size at most $k-1$ does not disconnect the graph. The maximum value of $k$ for which a graph $G$ is $k$-connected, is called the \emph{connectivity} of $G$ and is denoted by $\kappa(G)$.
We recall the statement of the well-known theorem by Menger next, and refer to the book by Diestel~\cite{DBLP:books/daglib/0030488} for this version of the theorem and its proof.
Note that the proof of Theorem~\ref{thm:menger} is constructive, meaning that the respective collection of paths can be computed in polynomial time~\cite{DBLP:books/daglib/0030488}.

\begin{theorem}[Menger's Theorem; Theorem 3.3.6,~\cite{DBLP:books/daglib/0030488}]\label{thm:menger}
A graph is $k$-connected if and only if it contains $k$ independent paths between any two vertices.
\end{theorem}

A graph is said to be \emph{bipartite} if its vertices can be partitioned into two sets $A$ and $B$ such that each edge in $G$ has its endpoints in different parts and we denote it as $G=(A\cup B, E)$.
A bipartite graph $G=(A\cup B, E)$ is called \emph{biconvex bipartite} if there exist  orderings $\sigma_A:A\rightarrow[|A|]$ and $\sigma_B:B\rightarrow[|B|]$ such that for each vertex $v\in A\cup B$, the vertices of $N(v)$ are consecutive under the corresponding ordering. These orderings together is called a \emph{biconvex ordering} of $A\cup B$. 
A bipartite graph $G=(A\cup B, E)$ is called \emph{convex bipartite} if there exists an ordering $\sigma:A\rightarrow[|A|]$ such that for every vertex $v\in B$, the vertices of $N(v)$ are consecutive under $\sigma$.
Such an ordering $\sigma$ is called \emph{convex ordering} of $A$ and can be computed in linear time~\cite{booth1976testing}.
For biconvex and convex bipartite graphs, given two vertices $u,v$ we denote by $u<v$ the fact that the ordering value of $u$ is less than the one of $v$ (similarly for $u\leq v$).
A superclass of convex bipartite graphs are \emph{chordal bipartite} graphs. 
We call a bipartite graph $G=(A\cup B, E)$, chordal bipartite, if it is $C_\ell$-free for $\ell\geq 6$.
A graph $G=(V,E)$ is called \emph{interval} if there is a set of intervals in the real line such that each vertex denotes an interval and two vertices are adjacent if and only if two intervals overlap.


\begin{definition}[\emph{Path decomposition}]
A \emph{path decomposition} of a graph $G = (V, E)$ is a sequence of subsets $B_1, \dots, B_\ell$ of $V$ (also called bags) satisfying the following conditions:
\begin{enumerate}
    \item $\bigcup_{i=1}^{\ell} B_i = V$.
    \item For every edge $uv \in E$, there exists an index $i$ such that $u, v \in B_i$.
    \item For every vertex $v \in V$, the set of indices $\{ i \mid v \in B_i \}$ forms a contiguous subsequence of $\{1, \dots, \ell\}$.
\end{enumerate}
The \emph{width} of a path decomposition is $\max_{1 \leq i \leq \ell} |B_i| - 1$, and the \emph{pathwidth} of $G$ is the minimum width over all possible path decompositions of $G$.
\end{definition}
\section{Missing Proofs from Section~\ref{section: All terminals on a tree}}
This section is dedicated to the proofs omitted from Section~\ref{section: All terminals on a tree}.
We restate the lemmata for convenience.

\lemmaAddTrees*
\begin{proof}
    Note that by considering iteratively a vertex in $N_{T_1}(\cup_{i=1}^kV_i)$ and adding it to an adjacent non full set $V_i$ the connectivity of each set is preserved.
    Moreover either one of the sets meets its size demand or $T_1\subseteq\cup_{i=1}^k V_i$. 
    If one of the sets meets its size demand, say $V_i$ then our algorithm outputs $V_i$ and an instance of the original problem for $k-1$, $G-V_i, C\setminus c_i, \mathcal T\setminus T_1$ since we have found one full set that intersects exactly one dominating tree. The graph created after removing this set is $(k-1)$-connected since it contains $k-1$ dominating trees. 

    Then since $T_1$ is dominating to $V\setminus \mathcal T$ and all of its vertices are already added to some vertex set we are able to add the vertices of $V\setminus \mathcal T$ one by one to an adjacent non full set again until either one of them becomes full in which case we restart the algorithm or until each of these vertices is added to some vertex set.
\end{proof}

\lemmaCorrectness*
\begin{proof}
    Since we have already argued for the polynomial running time of each algorithm called as a subroutine in Algorithm~\ref{case-alg:dominating trees}, its polynomial running time is also immediate.
    
    Notice also that during each call of another algorithm if we have found a smaller instance of the general problem then we stop our algorithm and output this and our lemma is satisfied.

    If this case does not occur and all the steps of Algorithm~\ref{case-alg:dominating trees} are performed.
    As argued in lemma~\ref{case-lemma: AddTrees} after the first step the only vertices that do not belong in any vertex set are the ones in $T_2,\ldots, T_k$.
    After the second step, due to lemma~\ref{case-lemma:ensure adjacencies in trees} each vertex set is either adjacent to enough available vertices to meet its size demand, or it has been assigned to a distinct tree and also contains a vertex from it.
    Then for the ones that have been assigned to some tree, using the vertex of that tree they already contain to ensure that there are more adjacent vertices of that tree in the third step algorithm Add Vertices either adds the whole tree to this vertex set or meets its size demand as proven in lemma~\ref{case-lemma: Sets in Under contain a tree}. 
    For the ones that have not been assigned some tree, there are enough available vertices assigned to each of them, so by adding them one by one each set becomes full.
    Hence lastly each non full set that remains contains a dominating tree and hence we can add available vertices to it arbitrarily until its full.
    Since every set is full and we have initialized each of them to contain the desired terminal vertex the partition our algorithm outputs is the desired one and satisfies our lemma for $l=k$.
    
\end{proof}

\section{Missing Part of Section~\ref{section: general case}}
\nameref{Description: Add vertices} is quite straightforward, both when it comes to its correctness and to the fact that it requires polynomial time so due to space restrictions we omit the corresponding pseudocode and provide directly the lemma that proves its correctness.
We refer to the algorithm realizing this step as $\text{Add Tree Vertices}(\mathcal T_0, \mathcal T_1, \mathcal T_{>1},V_1,\ldots, V_k,C, N, V(G)\setminus V(\mathcal T))$.


The two lemmata below prove the correctness of Algorithm Add Tree Vertices while its polynomial running time is easy to see as each vertex of $G$ is considered at most once to be added to some set.

\begin{restatable}{lemma}{lemmaVerticesOfTrees}\label{lemma: Add Vertices Of Trees} 
    Given $k$ dominating trees of a graph $G$ separated in three categories $\mathcal T_0, \mathcal T_1, \mathcal T_{>1}$, and $k$ distinct connected vertex sets, each containing one terminal vertex, Algorithm Add Tree Vertices adds each vertex of the trees in $\mathcal T_1\cup \mathcal T_{>1}$ to some vertex set such that each vertex set induces a connected graph. When this algorithm runs as part of Algorithm~\ref{alg:dominating trees-general} (\ref{Description: Add vertices}), then either each vertex of each tree in $\mathcal T_1\cup \mathcal T_{>1}$ is added to some vertex set inducing a connected graph or the algorithm restarts for a reduced value of $k$.
\end{restatable}
\begin{proof}
    When running Algorithm Add Tree Vertices as part of Algorithm~\ref{alg:dominating trees-general}, each tree that contains at least one terminal is considered separately.
    Then, each vertex of the tree that is adjacent to a non-full set is added to that vertex set, hence preserving the connectivity of each set formed.
    This procedure terminates when all the vertices of the considered tree are added to some vertex set, and each of those sets is not full or because one of those met its size demand.

    If it terminates due to the second condition, meaning some set $V$ becomes full, notice that $V$ is connected, and its removal reduces the connectivity of the graph by at most one because it does not contain any vertex from $k-1$ dominating trees, and hence, those trees cause the graph $G-V$ to be $(k-1)$-connected.
\end{proof}

\begin{restatable}{lemma}{lemmaUnassignedVertices}\label{lemma:Add Non Tree Vertices}
     Given $k$ connected vertex sets $\{V_1,\ldots, V_k\}$ such that their union induces a collection of connected dominating to $U$ sets and $k$ natural numbers $\{n_1,\ldots, n_k\}$, Algorithm Add Tree Vertices adds each vertex of $V(G)\setminus V(\mathcal T)$ to some neighboring set $V_i\in\{V_1,\ldots, V_k\}$. When this algorithm runs as part of Algorithm~\ref{alg:dominating trees-general} (\ref{Description: Add vertices}), then either each vertex of $V(G)\setminus V(\mathcal T)$ is added to some set inducing a connected graph or the algorithm restarts for a reduced value of $k$.
\end{restatable}
\begin{proof}
    The vertices of $V(G)\setminus V(\mathcal T)$ are considered one by one.
    Any such vertex is adjacent to some vertex set because, due to lemma~\ref{lemma: Add Vertices Of Trees}, the union of the vertex sets contains a dominating to $V(G)\setminus V(\mathcal T)$ set. 
    Each vertex of $V(G)\setminus V(\mathcal T)$ is added to a set, increasing its size by one and preserving its connectivity.
    This terminates when all vertices of $V(G)\setminus V(\mathcal T)$ are added to some connected vertex set or if one set meets its size demand.
    In the second case Algorithm~\ref{alg:dominating trees-general} restarts for a reduced value of $k$.

    If the loop terminates due to the second condition when run as part of Algorithm~\ref{alg:dominating trees-general}, the full set causing this termination, is connected. Also, its removal reduces the connectivity of the graph by at most one because it does not contain any vertex from $k-1$ dominating trees, and hence those trees make the graph $(G\setminus V_s)$ to be $(k-1)$-connected.
\end{proof}

In the next step we essentially find an $l$-connected subgraph of $G$, $G'$, $l$ dominating trees of $G'$ such that one of them contains $l$ terminal vertices such that the size demands of the corresponding vertex sets sum up to at most the size of $G'$. In other words in this step we find an instance of the special case described in the previous section in order to give it as an input to the corresponding algorithm afterwards.
This is done in polynomial time as each tree is considered at most once.
In lemma~\ref{lemma: Separate trees to sets} we prove the existence of such a subgraph.
\begin{algorithm}
\caption{{Choose a Tree Set (\ref{Description: Chose a tree set})}
    \label{alg:Separate trees to sets}}
    \begin{algorithmic}[1]
    \Require{A graph $G$, three sets of dominating trees $\mathcal T_0, \mathcal T_1, \mathcal T_{>1}$, a set of $k$ vertices $C=\{c_1,\ldots, c_k\}$, a set of $k$ natural numbers $N=\{n_1,\ldots, n_k\}$ and $k$ connected vertex sets $\{V_1,\ldots, V_k\}$ such that $c_i\in V_i$ for all $i\in [k]$}
    \While{$\mathcal T_0\neq\emptyset$}\label{Step: Separate trees to sets}
      \State Let $T\in \mathcal T_{>1}$ that contains $l$ terminals $\{c_1,\ldots, c_l\}\in C$ and let $V_1,\ldots, V_l$ contains $c_1,\ldots, c_l$, respectively.
      \State Pick $T_1',\ldots, T_{l-1}'\in \mathcal T_0$.
      \State $\mathcal T_{>1}\gets \mathcal T_{>1}\setminus T$
      \State $\mathcal T_0\gets \mathcal T_0\setminus\{T_1',\ldots, T_{l-1}'\}$.
      \If{$|T|+\sum_{i=1}^{l-1}|T_i'|\geq \sum_{i=1}^l n_i-|V_i|$}\label{Step: Condition for convenient tree set}
      \State \Return $T, T'_1,\ldots, T'_{l-1}$
      \EndIf
      \EndWhile
    \end{algorithmic}
\end{algorithm}
\begin{lemma}\label{lemma: Separate trees to sets}
    Let $G$ be a graph, $n_1,\ldots, n_k$ be  $k$ natural numbers such that $\sum_{i=1}^kn_i=|V(G)|$, $k$ dominating trees partitioned into three sets $\mathcal T_0,\mathcal T_1, \mathcal T_{>1}$, and $V_i$ be a connected vertex set containing terminal $c_i$, where $i\in [k]$, are given as an input to Algorithm~\ref{alg:Separate trees to sets} (\ref{Description: Chose a tree set}). Then, there exists a tree $T\in \mathcal{T}_{>1}$ with terminals $c_1,c_2,\ldots,c_l$ and $l-1$ trees $ T'_1,\ldots, T'_{l-1}\in \mathcal{T}_0$ such that $|T| + \sum_{i=1}^{l-1} |T_i'|\geq \sum_{i=1}^{l} n_i -|V_i|$.
\end{lemma} 
\begin{proof}
    Notice that when Algorithm~\ref{alg:Separate trees to sets} runs as a part of Algorithm~\ref{alg:dominating trees-general}, it essentially groups each tree containing more than one terminal with the respective number of trees containing no terminals such that the tree set contains an equal number of terminal vertices to dominating trees.
    We do this by removing $l-1$ trees from $\mathcal T_0$ any time a tree  $T \in \mathcal T_{>1}$ is considered, where $l$ is the number of terminal vertices in $T$.
    We know that enough such trees exist because we have a total of $k$ trees and $k$ terminals.

    Let $\mathcal P_1,\ldots, \mathcal P_s$ be those sets of trees (the ones originally in $\mathcal T_0$) and assume for a contradiction that none of those sets satisfies the condition of line~\ref{Step: Condition for convenient tree set}. 
    Moreover, due to Lemma~\ref{lemma:Add Non Tree Vertices}, we know that all the vertices not originally in $\mathcal T$ and the ones in some tree in $\mathcal T_1\cup \mathcal T_{>1}$ are contained in some non-full set.
     Denote the set of vertex sets created so far, meaning $\{V_1,\ldots, V_k\}$ by $\mathcal V$ and by $n_{\mathcal V}$ the sum of their size demands.
    Hence, $|V(G)|=\sum_{i=1}^s\sum_{P\in \mathcal P_i}|V(P)|+\sum_{H\in\mathcal V}|V(H)|$.
    Due to the previous observation and our assumption $|V(G)|<|V(G)|-n_{\mathcal V}+n_{\mathcal V}$, which leads us to the desired contradiction.
\end{proof}
We are now ready to put everything together to create the main algorithm of this section and prove its correctness.
\begin{algorithm}
\caption{{Gy\H{o}ri--Lov\'{a}sz theorem for graphs with $k$ dominating trees}
    \label{alg:dominating trees-general}}
    \begin{algorithmic}[1]
    \Require{A graph $G$, three sets of dominating trees $\mathcal T_0, \mathcal T_1, \mathcal T_{>1}$, a set of $k$ vertices $C=\{c_1,\ldots, c_k\}$, a set of $k$ natural numbers $N=\{n_1,\ldots, n_k\}$ and $k$ connected vertex sets $\{V_1,\ldots, V_k\}$ such that $c_i\in V_i$ for all $i\in [k]$}
    \State $(V_1,\ldots, V_k)\gets \text{Add Vertices}(\mathcal T_0,\mathcal T_1,\mathcal T_{>1}, V_1\ldots, V_k, C, N, U)$
    \State $ (T,T_1,\ldots, T_{l-1})\gets\text{Choose a Tree Set}(\mathcal T_0,\mathcal T_1,\mathcal T_{>1}, C, N, V_1,\ldots, V_k)$
    \State Algorithm~\ref{case-alg:dominating trees}$(\cup_{i=1}^lV_i\cup_{i=1}^{l-1}T'_i\cup T, T, T'_1,\ldots, T'_{l-1}, c_1,\ldots, c_l, n_1,\ldots, n_l)$.
    \end{algorithmic}
\end{algorithm}
Due to lemma~\ref{case-lemma: algorithm correctness} we conclude that Algorithm~\ref{alg:dominating trees-general} runs in polynomial time as we have already also argued that the first two steps of the algorithm also run in polynomial time.
The following lemma also implies theorem~\ref{alg:dominating trees-general} which concludes this section.
\begin{restatable}{lemma}{lemmaCorrectnessGeneral}
    Let a $k$-connected graph $G$, $k$ distinct terminal vertices $c_1,\ldots, c_k$, $k$ vertex disjoint dominating trees $T_1,\ldots, T_k$ and $k$ natural numbers such that $\sum_{i=1}^kn_i=|V(G)|$ are given as an input to Algorithm~\ref{alg:dominating trees-general}. Then Algorithm~\ref{alg:dominating trees-general} returns a partition of $G$ into $k$ sets $V_1,\ldots, V_k$ such that $|V_i|=n_i$ and $G[V_i]$ is connected for all $i\in[k]$ in polynomial time.
\end{restatable}
\begin{proof}
    Since we have already argued that Algorithm~\ref{alg:dominating trees-general} runs in polynomial time it suffices to show that it returns the desired partition.

    Notice that the graph $G'$ given to Algorithm~\ref{case-alg:dominating trees} as input is indeed $l$ connected for some $l\in[k]$ since it contains $l$ dominating trees, say $T_1,\ldots, T_l$. 
    Moreover, one of those trees, say $T_1$, contains $l$ terminals and the other none due to the way we group the trees together in algorithm Choose a Tree Set (one tree from $T_{>1}$ that contains $s$ terminal vertices and $s-1$ trees from $T_0$ ).
    
    Since the tree set chosen satisfies $\sum_{i=1}^l T_i\geq \sum_{i=1}^ln_i-|V_i|$ we get $|V(G')|\geq \sum_{i=1}^l n_i$. If $|V(G')|=\sum_{i=1}^l n_i$ it is immediate that we have the desired input so the algorithm outputs $l$ full sets and iterates on the smaller instance $G\setminus \cup_{i=1}^lV_i$, $\mathcal T\setminus\{T_1,\ldots, T_l\}$, $C\setminus\{c_1,\ldots, c_l\}$ and $n_{k-l},\ldots, n_k$.
    In the case where  $|V(G')|> \sum_{i=1}^l n_i$ we briefly adjust the size demands. Specifically we chose a vertex set, say $V_1$ and increase its demand from $n_1$ to $|V(G')|-\sum_{i=2}^ln_i$. After computing the $l$ partition through Algorithm~\ref{alg:dominating trees-general} we remove $|V(G')|-\sum_{i=1}^ln_i$ vertices from $V_1$ one by one ensuring that each time we do not destroy the connectivity and also do not remove $c_1$.
    Then we again output the $l$ sets and restart the algorithm for $G\setminus \cup_{i=1}^lV_i$, $\mathcal T\setminus\{T_1,\ldots, T_l\}$, $C\setminus\{c_1,\ldots, c_l\}$ and $n_{k-l},\ldots, n_k$.
    
\end{proof}
\def\P{\mathcal{P}}

\section{Gy\H{o}ri--Lov\'{a}sz Partition for Interval Graphs}\label{sec:interval}
In this section, we show that every $k$-connected interval graph admits a CDS partition of size $k$ (Theorem~\ref{thm:cds_interval}), which leads to an efficient construction for the Gy\H{o}ri--Lov\'{a}sz theorem on this graph class (Theorem~\ref{thm:convex_algo}). We now restate the main result of this section for convenience.

\cdsinterval*





\begin{algorithm}
\caption{{ $k$-CDS Interval Graph}
    \label{alg: k-cDSintervalGraph}}
    \begin{algorithmic}[1]
     \Require{A $k$-connected interval graph $G$} 
     \State Compute a minimum width path decomposition $B_1, \dots, B_\ell$ in $G$.
     \State $V(G)\gets V(G)\cup\{s,t\}$, $E(G)\gets E(G)\cup \{sv:\text{ for all } v\in B_1\}\cup\{tv:\text{ for all } v\in B_\ell\}$
     \State Compute $k$ disjoint $s-t$ paths in $G$, $\P = P_1, \dots, P_k$ and delete the endpoints $s$ and $t$ from each path.
     
     \State\Return $\P$
    \end{algorithmic}
    \end{algorithm}
%
    The algorithm takes advantage of a path decomposition\footnote{For the definition of path decomposition, refer to Appendix~\ref{sec:graph_notation}.} of the graph. To prove \cref{thm:cds_interval}, we show that the Algorithm~\ref{alg: k-cDSintervalGraph} produces $k$ dominating paths in polynomial time.

\begin{lemma}
   Algorithm~\ref{alg: k-cDSintervalGraph} returns $k$ disjoint paths, each dominating $V(G)$, and runs in polynomial time.
\end{lemma}
\begin{proof}
    It is well-known that in a minimum width path decomposition of an interval graph $G$, the induced subgraph of each bag $B_i$ forms a clique in $G$~\cite{DBLP:conf/wg/Garbe94}.
    Furthermore, each bag acts as a separator in $G$, meaning that a path from $s$ to $t$ must cross every bag. Combining both facts, we conclude that a path $P$ from $s$ to $t$, excluding the endpoints $s$ and $t$, dominates $V$.

    Since every bag $B_i$ is a separator, and by the $k$-connectivity of $G$, each bag has size at least $k$. Adding a vertex to a $k$-connected graph and connecting it to at least $k$ vertices that preserves $k$-connectivity. As a result, there are $k$ disjoint paths between $s$ and $t$.

  Computing a minimum-width path decomposition for interval graphs can be done in polynomial time~\cite{DBLP:conf/wg/Garbe94}. Computing $k$ disjoint paths can also be done in polynomial time by using a maximum flow algorithm. Thus, the algorithm runs in polynomial time.
\end{proof}
\def\T{\mathcal{T}}

\section{Gy\H{o}ri--Lov\'{a}sz Partition for Biconvex Graphs}\label{sec:biconvex}

This section is dedicated to the result for biconvex bipartite graphs. We show that every $k$-connected biconvex bipartite graph admits a CDS partition of size $k$ (Theorem~\ref{thm:cds_biconvex}), which leads to an efficient construction for the Gy\H{o}ri--Lov\'{a}sz theorem on this graph class (\cref{thm:biconvex_algo}). We now restate the main result of this section for convenience.

\cdsbiconvex*



\begin{algorithm}
\caption{{$k$-CDS Biconvex}
    \label{alg: k-cDSbiConvex}}
    \begin{algorithmic}[1]
     \Require{A $k$-connected biconvex graph $G = (A \cup B, E)$ with biconvex ordering $\sigma$ such that $\sigma_A = (a_1, \dots, a_n)$ and $\sigma_B =(b_1, \dots, b_m)$} 
     \State  Compute $k$ disjoint $a_1-a_n$ induced paths $\T = T_1, \dots, T_k$ and delete the endpoints $a_1$ and $a_n$ from each path \label{step::biconvexInit}
     \For{every $T_i\in\T$}\label{step::biconvexMakeItDS}
     \If{$T_i\cap N(b_1)=\emptyset$}
     \State $T_i\gets T_i\cup v$ such that $v\in N(b_1)$ and $v\not\in V(\T)$
     \EndIf
     \If{$T_i\cap N(b_m)=\emptyset$}
     \State $T_i\gets T_i\cup v$ such that $v\in N(b_m)$ and $v\not\in V(\T)$
     \EndIf
     \EndFor
     \State\Return $\T$
    \end{algorithmic}
    \end{algorithm}

To prove the Theorem~\ref{thm:cds_biconvex}, we design an algorithm to construct a CDS partition for biconvex graphs, given in Algorithm~\ref{alg: k-cDSbiConvex}. To prove the correctness of Algorithm~\ref{alg: k-cDSbiConvex}, first we need to show that there always exists a neighbor of $b_1$ or $b_m$ in $V(G) \setminus V(\T)$ to add, at the lines 4 and 6.

\begin{restatable}{lemma}{lemmaEnsureNeighborsofB}
    \label{lemma::atMostOneNeighbor}
    Let $G$ be a connected biconvex graph $G = (A \cup B, E)$ and $\sigma$ be a biconvex ordering such that $\sigma_A = (a_1, \dots, a_n)$ and $\sigma_B = (b_1, \dots, b_m)$.  
    Let $P'$ be an induced $(a_1,a_n)$-path and $P=P'\setminus \{a_1,a_n\}$. Then the path $P$ contains at most one neighbor of $b_1$ and at most one neighbor of $b_m$. Furthermore, the path $P$ dominates the vertex set $A$.
\end{restatable}
\begin{proof}
    Let $P$ be such a path and assume that $|V(P) \cap N(b_1)| > 1$.  
    Let $a'_1, a'_2 \in N(b_1)\cap V(P)$.  
    Since $G$ is bipartite, $a'_1, a'_2 \in A$.  
    Let $b'_1, b'_2, b'_3, b'_4 \in B$ be the neighbors of $a'_1, a'_2$ in $P$, such that $b'_1 \leq b'_2 \leq b'_3 \leq b'_4$.  
    Note that the $b'_i$'s are not necessarily distinct; however, at most one of them has $a'_1$ and $a'_2$ as neighbors in $P$.  
    Without loss of generality, let $a'_1$ be the neighbor of $b'_4 = b_x$.  
    Consequently, the neighborhood of $a'_1$ in $G$ contains $b_1, \dots, b_x$ by convexity, and in particular, the vertices $b'_1, b'_2, b'_3, b'_4$.  
    Thus, the degree of $a'_1$ is at least three in $G[V(P)]$, contradicting the precondition that $P$ is an induced path.  
    
    Symmetrically, we can make the same argument for $b_m$.  
    Suppose $|V(P) \cap N(b_m)| > 1$.  
    Let $a''_1, a''_2 \in N(b_m)$ be two vertices in $P$.  
    Since $G$ is bipartite, $a''_1, a''_2 \in A$.  
    Let $b''_1, b''_2, b''_3, b''_4 \in B$ be the neighbors of $a''_1, a''_2$ in $P$, such that $b''_1 \leq b''_2 \leq b''_3 \leq b''_4$.  
    Note that the $b''_i$'s are not necessarily distinct; however, at most one of them has $a''_1$ and $a''_2$ as neighbors in $P$.  
    Without loss of generality, let $a''_1$  be the neighbor of $b''_1 = b_y$.  
    Consequently, the neighborhood of $a''_1$ in $G$ contains $b_y, \dots, b_m$ by convexity, and in particular, the vertices $b''_1, b''_2, b''_3, b''_4$.  
    Thus, the degree of $a''_1$ is at least three in $G[V(P)]$, contradicting the precondition that $P$ is an induced path.

    Next, we want to show that the path $P$ dominates the vertex set $A$. Clearly, $a_1$ and $a_n$ are dominated by $P$.  
    Assume there exists a vertex $a_i \in A \setminus \{a_1, a_n\}$ that is not dominated by $P$.  
    Originally, the path extended from $a_1$ to $a_n$, which implies the existence of consecutive vertices $a_x, b_y, a_z$ in the original $P$ with $a_x < a_z$ such that $a_x \leq a_i \leq a_z$.  
    However, by convexity $a_i \in N(b_y)$, contradicting our assumption.
\end{proof}

We now use the fact that $k$-connectivity ensures each vertex has at least $k$ neighbors.  
By \cref{lemma::atMostOneNeighbor}, each $T \in \T$ contains at most one neighbor of $b_1$ or $b_m$.  
Consequently, if $T$ has no neighbor of $b_1$ or $b_m$, there must be an available neighbor in $A$ in step~\ref{step::biconvexMakeItDS} as $|\T| = k$.



The sets in $\mathcal{T}$ are disjoint by construction. 
Recall that every $k$-connected graph guarantees the existence of at least $k$ disjoint paths between any pair of vertices.
Since every path $T$ computed in step~\ref{step::biconvexInit} is connected and already dominates $A$ (see \cref{lemma::atMostOneNeighbor}), every additional vertex in $A$ that is added to $T$ preserves its connectivity.
Next, we show that each set in $\mathcal{T}$ dominates the whole graph.

\begin{lemma}
   Let $\mathcal{T}$ be a returned set by Algorithm~\ref{alg: k-cDSbiConvex}. Then, every vertex set $T\in \mathcal{T}$ dominates the input graph $G$.
\end{lemma}
\begin{proof}
    Let $T \in \T$. By \cref{lemma::atMostOneNeighbor}, $T$ dominates $A$ as $T$ contains an $(a_1,a_n)$-path except $\{a_1,a_n\}$. Now, it remains to show that $T$ dominates $B$. 
    Let $b_x$ and $b_y$ be the vertices with the lowest and highest ordering value in $V(T) \cap B$.  
    First, consider the case $b_x = b_y$.  
    By the construction of $T$ in the algorithm, $N(V(T))$ contains $b_1$ and $b_m$, and by \cref{lemma::atMostOneNeighbor}, $b_x$ dominates $A$.  
    Consequently, $T$ dominates $B$, since there exists a path $b_1, a, b_x$ and a path $b_x, a', b_m$, where $a, a' \in V(T)$.
    
    Now, consider the case $b_x < b_y$.  
    We begin by proving that the sequence $b_x, \dots, b_y$ is dominated by $T$.  
    Note that $b_x$ and $b_y$ are in $T$ before loop at line~\ref{step::biconvexMakeItDS}, since in this step only vertices from $A$ are added.  
    That is, $b_x$ and $b_y$ are part of the path computed in step~\ref{step::biconvexInit}, and therefore, there exists a path $P$ in $T$ from $b_x$ to $b_y$.  
    Suppose there exists a vertex $b_i \in \{b_x, \dots, b_y\}$ that is not dominated by $P$.  
    Since there is a path from $b_x$ to $b_y$, there exist consecutive vertices $b_j, a_z, b_k$ in $P$ with $b_j < b_k$ such that $b_j \leq b_i \leq b_k$.  
    However, by convexity $b_i \in N(a_z)$, contradicting our assumption.  
    Finally, by construction $N(V(T))$ contains $b_1$ and $b_m$, and by \cref{lemma::atMostOneNeighbor}, $T$ dominates $A$.  
    Consequently, $T$ dominates $B$, since there exists a path $b_1, a, b$ and a path $b', a', b_m$, where $a, a',b,b' \in V(T)$.
\end{proof}

   By using a maximum flow algorithm, we can compute $k$ disjoint paths between $a_1$ and $a_n$ in polynomial time. In the final step, we can check if a $T \in \T$ has no neighbor in $b_1$ or $b_m$ and add extra vertex in $T$, in polynomial time. Hence,  Algorithm~\ref{alg: k-cDSbiConvex} runs in polynomial time. The collection returned by Algorithm~\ref{alg: k-cDSbiConvex}, consists of $k$ disjoint connected dominating sets, which can be extended to a CDS partition of $G$. Hence, \cref{thm:cds_biconvex} follows.
\section{4-Approximation for Convex Bipartite Graphs}\label{sec:convex-bip_partition}

This section is dedicated to the results for convex bipartite graphs. We show that every $4k$-connected convex bipartite graph admits a CDS partition of size $k$ (Theorem~\ref{thm:cds_convex}), which leads to the $4$-approximate efficient construction for the Gy\H{o}ri--Lov\'{a}sz theorem on this graph class (Theorem~\ref{thm:convex_algo}). We now restate the main result of this section for convenience.

\cdsconvex*

The remainder of this section is dedicated to the proof of Theorem~\ref{thm:cds_convex}. Intuitively, the proof strategy is as follows. We would like to construct each of the connected dominating sets based on a path reaching from the smallest to the largest vertex, in terms of the convex ordering. We can get $k$ such paths directly from Menger's theorem; by convexity, every such path automatically dominates the part $A$. It does not necessarily dominate the whole of $B$ however, and we will need to attach selected vertices of $A$ to the path, in order to achieve domination on $B$. To make this possible, we need that the paths traverse the graph ``fast enough'', so that enough neighbors are left for each vertex of $B$ to be dominated by each of the constructed connected dominating sets. We show that making the paths induced is sufficient to achieve a gap of order $k$, on average, between consecutive vertices of $A$ on each path. The final part of the proof is then to argue that the gaps left by the paths are wide enough so that the remaining vertices of $A$ can be assigned to different sets, with every vertex of $B$ being dominated by each of the sets.

We now proceed with the formal proof. 
Let $G = (A \cup B, E)$ be the input graph and $\sigma$ be the convex ordering defined on $A$ which can be computed in linear time~\cite{glover1967maximum}.
Let $|A| = \ell$, and let $A = \{a_1, \ldots, a_\ell\}$ be the indexing of the elements of $A$ that agrees with the convex ordering $\sigma$; that is, for each $1 \le i < j \le \ell$, $\sigma(a_i) < \sigma(a_j)$. Let $\mathcal{P} = \{P_1, \ldots, P_k\}$ be the $k$ paths constructed by Theorem~\ref{thm:menger}, invoked on $G$ with the endpoints $a_1$ and $a_\ell$. 
Additionally, we assume that the paths are induced. This is safe since if any path in $\mathcal{P}$ is not induced, it can be shortened via the edge violating the induced property, until there are no more such edges.

\begin{claim}\label{claim:domA}
    $\mathcal{P}$ is a collection of internally-vertex-disjoint induced paths in $G$. Additionally, for each path $P \in \mathcal{P}$, each vertex of $A$ is adjacent to some vertex of $P$.
\end{claim}
\begin{proof}
    The first part of the statement follows by construction. For the second part, consider a path $P \in \mathcal{P}$, and let $b_1, \ldots b_p$ be the sequence of vertices of $B$ appearing on $P$, in the order they appear on the path. By construction, $b_1$ is adjacent to $a_1$, and $b_p$ is adjacent to $a_\ell$. We now claim that $\bigcup_{i = 1}^p N(b_i) = A$. Indeed, by convexity, each $N(b_i)$ is a continuous interval of $A$. Since the union of intervals in an interval, and both the leftmost and the rightmost elements of $A$ are contained in $\bigcup_{i = 1}^p N(b_i)$ from $a_1 \in N(b_1)$, $a_\ell \in N(b_p)$, the union must contain the whole of $A$ as well.
\end{proof}

Next, we show that the fact that the paths in $\mathcal{P}$ contain no ``shortcuts'', together with the convex property, implies that the arrangement of each paths' vertices on $A$ satisfies some useful properties. 

\begin{claim}\label{claim:gap}
    Let $P$ be a path in $\mathcal{P}$, and let $a_{i_1} = a_1$, $a_{i_2}$, \ldots, $a_{i_q} = a_\ell$ be the sequence of vertices of $A$ appearing on $P$, in the order they appear on the path. The following two properties hold.
    \begin{enumerate}
        \item The vertices of $A$ on $P$ are monotone, i.e., $i_1 < i_2 < \ldots < i_q$.
        \item Any interval of $4k$ consecutive vertices of $A$ contains at most $3$ vertices of $P$.
    \end{enumerate}
\end{claim}
\begin{proof}
    We start by showing the first property. Assume that it does not hold, and consider the smallest index $j$ so that $i_{j + 1}$ violates the property. That is, $j \in [q - 1]$ is such that $i_1 < \ldots < i_j$, but $i_{j + 1} < i_j$. Since $i_1 = 1$, $j$ is at least $2$. We now reach the contradiction by considering two cases depending on the position of $i_{j + 1}$ compared to $i_{j - 1}$. If $i_{j - 1} < i_{j + 1} < i_j$, consider the vertex $b \in B$ that lies on $P$ between $a_{i_{j - 1}}$ and $a_{i_j}$. Since $b$ is adjacent to each of them, it is also adjacent to $a_{i_{j + 1}}$ by convexity, which contradicts the fact that $P$ is an induced path. Otherwise, $i_{j + 1} < i_{j - 1} < i_j$, and consider the vertex $b' \in B$ that lies on $P$ between $a_{i_j}$ and $a_{i_{j + 1}}$. Again, by convexity $b'$ is adjacent to $a_{i_{j - 1}}$, which is a contradiction, since $P$ is an induced path.

    We move to the second property. Specifically, we show that for each $j \in [q - 3]$, $i_{j + 3} - i_j > 4k$. Clearly, this implies the desired statement, since any interval of vertices of $A$ of length $4k$ containing at least $4$ vertices of $P$, contains also $4$ consecutive vertices of $P$ in $A$, violating the former property. Assume now that there exists $j \in [q - 3]$ with $i_{j + 3} - i_j \le 4k$. Let $b \in B$ be the vertex of $P$ between $a_{i_{j + 1}}$ and $a_{i_{j + 2}}$. Since $P$ is an induced path, $b$ is not adjacent to $a_{i_j}$ and $a_{i_{j + 3}}$.
    On the other hand, since $G$ is $4k$-connected, $b$ has at least $4k$ neighbors, which must form a consecutive interval of $A$. This interval includes $a_{i_{j + 1}}$ and is therefore between $a_{i_{j}}$ and $a_{i_{j + 3}}$, excluding these two vertices. This is a contradiction, since  $i_{j + 3} - i_j \le 4k$, so there are at most $4k - 1$ vertices between $a_{i_{j}}$ and $a_{i_{j + 3}}$ in $A$.
\end{proof}

We now complete the construction by distributing the elements of $A$ not covered by the paths in $\mathcal{P}$ evenly between the dominating sets. More specifically, let $a_{i_1}, \ldots, a_{i_t}$ be the vertices of $A' = A \setminus \bigcup_{P \in \mathcal{P}} P$. We partition $A'$ into $k$ sets $S_1$, \ldots $S_k$ by assigning these vertices alternatingly in the order of $A$. That is, $S_1$ receives each $k$-th vertex of $A'$, $S_2$ each $k$-th vertex of $A$ starting from $a_{i_2}$, and so on; formally, $S_i = \{a_{i_j} : j \equiv i\  (\text{mod } k)\}$. We now argue that each set $S_i$ is dominating to the part $B$.

\begin{claim}\label{claim:domB}
    For each $i \in [k]$, $S_i$ dominates $B$.
\end{claim}
\begin{proof}
    Let $b$ be a vertex of $B$ and $i$ be an index in $[k]$. Since $G$ is $4k$-connected, $b$ has at least $4k$ neighbors.  By convexity, these neighbors form a consecutive interval in $A$.  Let $W$ be a subset of exactly $4k$ neighbors of $b$ that forms a consecutive interval in $A$ as well. By Claim~\ref{claim:gap}, each of the $k$ paths in $\mathcal{P}$ contains at most $3$ vertices of $W$, therefore at least $k$ vertices of $W$ are not contained in any of the paths in $\mathcal{P}$. Therefore, $W$ contains $k$ consecutive vertices of $A'$. By construction, one of them belongs to $S_i$, concluding the proof.
\end{proof}

We are now ready to conclude the proof of the theorem. For each $i \in [k]$, construct a set $T_i = V(P_i) \setminus \{a_1, a_n\} \cup S_i$. We claim that $\mathcal{T} = \{T_1, \ldots, T_k\}$ is a packing of connected dominating sets in $G$. Indeed, for each $i \in [k]$, $T_i$ is a dominating set in $G$: $V(P_i) \setminus \{a_1, a_n\}$ dominates $A$ by Claim~\ref{claim:domA}, and $S_i$ dominates $B$ by Claim~\ref{claim:domB}. Furthermore, $G[T_i]$ is connected. This holds since $V(P_i) \setminus \{a_1, a_n\}$ is a subpath of the path $P_i$, therefore this set induces a connected subgraph, and each vertex of $S_i$ is adjacent to $V(P_i) \setminus \{a_1, a_n\}$ by Claim~\ref{claim:domA}. All sets in $\mathcal{T}$ are vertex-disjoint, since the paths in $\mathcal{P}$ are internally-vertex-disjoint, and the sets $\{S_i\}_{i \in [k]}$ form a partition of $A \setminus \bigcup_{P \in \mathcal{P}} P$. Finally, $\mathcal{T}$ can be made into a CDS partition of $G$ of size $k$ by assigning the remaining vertices $a_1$, $a_n$ to $T_1$. It is also clear that $\mathcal{T}$ can be constructed in polynomial time.

\end{document}